%% file: confn.tex
\documentclass[11pt,a4paper]{amsart}
\usepackage{amssymb}
\usepackage{epsfig}
\begin{document}
\newtheorem{thrm}{Theorem}
\newtheorem{thmf}{Th{\'e}or{\`e}me}
\newtheorem{thmm}{Theorem}
\def\thethmm{\ref{twist}$'$}
\newenvironment{thint}[1]{{\flushleft\sc{Th{\'e}or{\`e}me}}
      {#1}. \it}{\medskip} 
\newenvironment{thrm.}{{\flushleft\bf{Theorem}}. \it}{\medskip} 
\newenvironment{thrme}[1]{{\flushleft\sc{Theorem}}
      {#1}. \it}{\medskip} 
\newenvironment{propint}[1]{{\flushleft\sc{Proposition}}
      {#1}. \it}{\medskip} 
\newenvironment{corint}[1]{{\flushleft\sc{Corrolaire}}
      {#1}. \it}{\medskip} 
\newtheorem{corr}{Corollary}
\newtheorem{corf}{Corollaire}
\newtheorem{prop}{Proposition}
\newtheorem{defi}{Definition}
\newtheorem{deff}{D{\'e}finition}
\newtheorem{lem}{Lemma}
\newtheorem{lemf}{Lemme}
\def\fy{\varphi}
\def\ul{\underline}
\def\obsf{{\flushleft\bf Remarque. }}
\def\obs{{\flushleft\bf Remark. }}
\def\R{\mathbb{R}}
\def\C{\mathbb{C}}
\def\Z{\mathbb{Z}}
\def\N{\mathbb{N}}
\def\P{\mathbb{P}}
\def\Q{\mathbb{Q}}
\def\p{\pi_1(X)}
\def\Re{\mathrm{Re}}
\def\Im{\mathrm{Im}}
\def\H{\mathbf{H}}
\def\Hs{{Nil}^3}
\def\L{L_\alpha}
\def\h{\frac{1}{2}}
\def\M{\P(E)\times\P(E)^*\smallsetminus\mathcal{F}}
\def\m{\cp2\times{\cp2}^*\smallsetminus\mathcal{F}}
\def\g{\mathfrak{g}}
\def\l{\mathcal{L}}
\def\k{\mathfrak{h}}
\def\s3{\mathfrak{s}^3}
\def\nb{\nabla}
\def\Re{\mathrm{Re}}
\def\Im{\mathrm{Im}}
\newcommand\cp[1]{\mathbb{CP}^{#1}}
\renewcommand\o[1]{\mathcal{O}({#1})}
\renewcommand\d[1]{\partial_{#1}}
\def\sq{\square_X}
\def\gp{\dot\gamma}
\def\j{\mathcal{J}}

\title[Null-geodesics and the LeBrun correspondence]{Null-geodesics in
  complex conformal manifolds and the LeBrun
  correspondence}  
\author{Florin Alexandru Belgun}
\thanks{AMS classification~: 53C21, 53A30, 32C10.}
\thanks{Partially supported by the SFB 288 of the DFG}
\date{February 9, 2000.}

\begin{abstract}
In the complex-Riemannian framework we show
that a conformal manifold containing a compact,
simply-connected, null-geodesic is conformally flat. In dimension 3 we
use the LeBrun correspondence, that views a conformal 3-manifold
as the conformal infinity of a
selfdual four-manifolds. We also
find a relation between the conformal invariants of the conformal
infinity and its ambient. 
\end{abstract}

\maketitle

\section{Introduction}
On a complex manifold, the existence of a {\it complex-Riemannian
  metric} implies, in general, strong topological
  assumptions, especially if the manifold is compact ({\it e.g.} the
  --- square of the --- canonical bundle has to be trivial). However,
  any analytic (pseudo-) Riemannian (or conformal) manifold can be complexified,
  and a natural question is to see to what extent the {\it global}
  properties of the real manifold ({\it e.g.} existence of closed
  (null-) geodesics) hold for the complexified spaces. This
  complexification procedure naturally occurs in {\it twistor theory}
  (see below), which has been intensively studied for Riemannian
  space-times (see, {\it e.g.}, \cite{ahs}, \cite{hi}, \cite{leb3},
  \cite{ward}); the complex-Riemannian setting, in which historically
  the twistor theory was first introduced \cite{p1}, can provide a
  link to the Lorentzian geometry. 

In {\it complex conformal geometry} (which implies weaker assumptions
on the topology of the manifold), the conformal structure is
determined by the set of {\it null-geodesics}, which can be organized
as a complex manifold under some topological conditions \cite{leb1},
\cite{leb5}. A natural question is which complex conformal manifolds
admit compact null-geodesics; for example, if a {\it self-dual} manifold admits a
globally-defined {\it twistor space}, then application of a twistorial
interpretation of the Weyl tensor \cite{FW}, implies that it is
conformally flat, and the compact null-geodesic is
simply-connected. 
\smallskip 
 
Our main result (section 4, Theorem
\ref{comp3}) states that, if a conformal complex $n$-manifold admits 
a rational curve as a null-geodesic, then it is conformally
flat (see also \cite{ye} for the case of a complex projective
manifold).  The proof uses
the properties of  Jacobi fields along the considered compact,
simply-connected, null-geodesic~: namely, we compute the normal 
bundle of a compact, simply-connected, null-geodesic, and we show that
the small deformations of the latter as a compact curve, or as a
null-geodesic, coincide (section 4, Proposition \ref{deform}).
In addition to that, we use, for the (more difficult)
case of dimension 3, a criterion for conformal flatness from
\cite{FW}, and we apply it to a {\it locally defined}, by the {\it LeBrun
correspondence} (see below), self-dual ambient.
\smallskip

The other topic of this paper uses implicitly another application of
twistor theory: It has been shown by LeBrun \cite{leb1}, \cite{leb2}, that any
conformal 3-manifold can be locally realized as the {\it conformal
  infinity} of a self-dual Einstein (with non-zero scalar curvature)
4-manifold. We have, thus, a local correspondence assigning to a conformal
structure in dimension 3 a self-dual Einstein metric in dimension
4, which we call the {\it LeBrun correspondence}. 

As conformal structures of both manifolds are encoded in the complex,
resp. $CR$, structure of their twistor spaces, they are implicitly
related, for example if the 3-manifold $M$ is conformally flat, its
ambient $N$ equally is. It is, however, difficult to obtain an
explicit relation between the conformal invariants of $M$ and those of
$N$ by twistorial methods, as there is no simple expression of the
{\it Cotton-York tensor} of $M^3$ in twistorial terms, and the
twistorial interpretation of the {\it Weyl tensor} of $N^4$ is highly
non-linear \cite{FW}.

In this paper we find a relation between these two
conformal invariants of the manifolds involved in the LeBrun
correspondence, or, more generally, of an umbilic submanifold $M^3$ and
of its self-dual ambient $N^4$. It appears that the Weyl tensor of $N^4$
identically vanishes along $M^3$, and thus the Cotton-York tensor of
$N^4$, restricted to $M^3$, is conformally invariant and can be
identified with the Cotton-York tensor of $M^3$; in this case, it is
also equal to the normal 
derivative of the Weyl tensor of $N^4$ (section 3, Theorem
\ref{umbilic}). This gives conditions for an open  self-dual 4-manifold to
admit a conformal infinity.
\bigskip




The paper is organized as follows~: in section 2 we recall a few basic
facts about complex- Riemannian and -conformal geometry, in section 3
we relate the conformal invariants of a 3-dimensional conformal
infinity to those of its self-dual ambient (arising from the LeBrun
correspondence), and in section 4 we state our results about conformal
complex manifolds containing compact null-geodesics.
\smallskip

Throughout the paper we use the following conventions: in
complex-Riemannian (or -conformal) geometry we use the same
terminology as in the real framework (metric, Levi-Civita connection,
curvature), 
and the holomorphic bundles are denoted like the corresponding bundles
in real geometry (for example, the holomorphic tangent bundle of {\bf
  M} is denoted simply by $T\mathbf{M}$, rather than the more precise
$T^{1,0}\mathbf{M}$)~; manifolds with holomorphic conformal structures
are denoted by bold-face letters (except in section 3, where the
results hold also in the real framework).


\section{Holomorphic conformal geometry}
\begin{defi} Let $\mathbf{M}$ be a complex manifold, let $n$ be its complex
  dimension.  A {\em complex-Riemannian metric} $g$ on it is a
    holomorphic section of $S^2T^*\mathbf{M}$ which is non-degenerate
    at any point. A {\em holomorphic conformal structure} on
    $\mathbf{M}$ is a holomorphic line subbundle $C$ in
    $S^2T^*\mathbf{M}$  such that any non-vanishing local section of
    $C$ is a local complex-Riemannian metric.
\end{defi}

During the rest of this section, and of the whole fourth section of
this paper, we shall simply denote these structures as {\it metric} and
{\it conformal structure} (therefore omitting any reference to the
complex framework).

A metric on $T\mathbf{M}$ induces metrics (i.e. non-degenerate
symmetric bilinear forms) on all tensor bundles, in particular the
square of the canonical bundle $\kappa:=\Lambda^nT^*\mathbf{M}$ is
trivialized. If we can choose {\it an orientation} (defined as
follows), then the canonical bundle itself can be trivialized by a {\it
volume form} of norm 1. There are exactly 2 such forms in each fiber
of $\Lambda^nT^*\mathbf{M}$, and an {\it orientation} is the choice of
one of them, depending continously (thus holomorphically) on the base
point.

\obs The notion of {\it orientation} is generally related to the
reduction (when possible) of the structure group of the frame bundle
from $G$ to the connected component of the identity $G_0$; in absence
of any structure, the group $G$ is simply the connected $GL(n,\C)$, so the
notion of orientation has no meaning in ``raw'' complex geometry. But
a Riemannian metric on 
$\mathbf{M}$ is equivalent to the reduction 
of its frame bundle to a $O(n,\C)$-bundle, where $O(n,\C):=\{A\in
GL(n,\C)|A^tA=\mathbf{1}\}$; a further choice of an
orientation reduces the structure group of the frame bundle to
the connected component of this group, containing the identity:
$SO(n,\C)~:=O(n,\C)\cap SL(n,\C)$.

Unlike in the real framework, these reductions, always possible
on (small) contractible open sets, are submitted to some topological
constraints if we want to define them globally on $\mathbf{M}$.
\smallskip

Weaker constraints are implied by the existence of a conformal
structure~: the square of the canonical bundle needs to admit a $n$th
order root $L^{-2}\simeq C$, where $L:=\kappa^{-1/n}$ is
the\footnote{the dual of $L$ is a square root of $C$; 
  the choice of such a square root is implied, if $n=2m+1$, by
  the conformal structure $C$ as $\kappa=C^m\otimes L^{-1}$; if $n$ is
  even, neither $C$ nor an orientation --- see below --- imply the
  choice of $L$.}
(respectively, one of the) {\it weight
  1--density bundle} of the manifold $\mathbf{M}$. We can study the conformal
structure $C$ using the formalism of density bundles and of Weyl
derivatives \cite{gconf}. From now on, we shall not make
use of the weight 1--density bundle, but only of $C=L^{-2}$, which is
enough to define the conformal structure.

\obs A conformal structure is equivalent to the reduction of the
structure group of the frame bundle to $CO(n,\C):=O(n,\C)\times
\C^*/\{\pm\mathbf{1}\}$ (the quotient is due to the fact that
$-\mathbf{1}\in O(n,\C)$). This group is disconnected if $n$ is even,
and the connected component of {\bf 1} is
$CO_0(n,\C)=SO(n,\C)\times\C^*/\{\pm\mathbf{1}\}$, but it is connected
if $n$ is odd (and in this case the right hand side of the previous
equality coincides with $CO(n,\C)$). Therefore, although a complex-Riemannian
3-manifold admits, locally, 2 possible orientations, they are
conformally equivalent, fact that makes impossible a canonical way
to associate an orientation to a metric in the conformal class.


\obs In complex-, as in real-Riemannian geometry, the orientation
determines (and is determined by) a family of compatible oriented orthonormal
basis in $T\mathbf{M}$; if $\dim\mathbf{M}$ is 
even, by multiplying all the
vectors of such a basis by a non-zero complex number we obtain an
oriented orthonormal basis compatible with another metric in the
conformal class
(if $\dim\mathbf{M}$ is odd, multiplication by $-1$ yields a basis
compatible with the same metric, but with the opposite orientation). 

The notion of orientation, in four-dimensional conformal geometry, is
important for the definition of {\it anti--}, resp. {\it self-duality}
(see below). More generally, the {\it Hodge $*$ operator} (defined as
in real Riemannian geometry) is conformally invariant, and gives an
explicit expression for the splitting of the bundle of 2-forms and of
the curvature tensor \cite{st}, see also next section.
\medskip

Maybe the simplest way to view a conformal structure is as an
equivalence class $c$ of {\it local} metrics, two such local representants
$g$ and $h$ satisfying, on the open set where both are defined, to
$g=fh$, with $f$ a non-vanishing holomorphic function. Unlike in the
real framework, global representants may not exist in general. From
now on we shall
consider a conformal structure on $\mathbf{M}$ as being given by the conformal
class $c$ rather than by the line bundle $C$.

Geometrically, a conformal structure is given by its {\it isotropy
  cone} $\mathcal{C}\subset T\mathbf{M}$ of vectors of norm 0. Because of the
  non-degeneracy of any local metric in $c$, the projective isotropy
  cone $\P(\mathcal{C})$ is a non-degenerate hyperquadric in $\P(T\mathbf{M})$. In
  dimension 3, $\P(\mathcal{C})$ is a conic (curve) in $\cp2$, and in dimension
  4 it is a ruled surface in $\cp3$~; therefore, in this latter case,
  there are 2 families 
  --- each of which can be characterized with respect to a given orientation
  \cite{FW} --- of isotropic planes in $T\mathbf{M}$ called $\alpha$-,
  resp. $\beta${\it -planes}.

For any local metric we have a {\it Levi-Civita connection}, whose
curvature has the same components as in the real Riemannian geometry
(see next section). In particular, the {\it Weyl tensor} is
independent of the metric in the conformal class.

The geodesics for which the tangent direction at a point (and thus, at
any point) is isotropic are called {\it null-geodesics}, and they are
locally independent (up to a reparametrization) of the metric in the
conformal class. The same is true for higher-dimensional totally
geodesic and isotropic submanifolds --- if they exist ---, called
null-submanifolds. In dimension 4 they are $\alpha$-,
resp. $\beta${\it -surfaces} (tangent to $\alpha$-, resp. $\beta$-planes,
see above), and they exist if and only if the (oriented) conformal
structure is {\it anti-}, resp. {\it self-dual}, i.e. the component
$W^+$, resp. $W^-$, of the Weyl tensor $W$ of $(\mathbf{M},c)$ vanishes
identically \cite{ahs},\cite{p1},\cite{FW}.

If the manifold is self-dual, one considers locally the twistor space
$Z$ of $(\mathbf{M}^4,c)$ as the set of $\beta$-surfaces\footnote{see \cite{leb1},
\cite{leb5}, \cite{FW} and section 4 for an explanation of the
difficulties of a global definition of the twistor space.}. It is a
3-dimensional complex manifold \cite{ahs},\cite{p1} containing
rational curves whose normal bundle is isomorphic to $\o1\oplus\o1$
(called {\it twistor lines}) (where $\o1$ is the dual of the {\it
  tautological bundle} $\o{-1}$ of $\cp1$). If, in addition, we can choose an
Einstein metric $g$ in the conformal class $c$, we get an extra structure
on $Z$, namely a distribution of 2-planes, which is totally integrable
(and yields a foliation) if the scalar curvature of $g$ vanishes,
otherwise it is a contact structure \cite{ward},\cite{hi},\cite{leb1}.
Conversely, from a manifold $Z$ containing twistor lines as above
(called a {\it twistor space}), plus ---
possibly --- the additional 2-planes distribution, one can recover
--- at least locally ---,
via the {\it reverse Penrose construction} \cite{ahs}, the self-dual
manifold $(\mathbf{M}^4,c)$.

In all generality, one can always consider locally (on a geodesically
convex open set, for example, see section 4) the space of
null-geodesics of a conformal manifold $(\mathbf{M}^n,c)$, 
and the key point in the {\it LeBrun correspondence} (defined below) is that the
space of null-geodesics of a 3-dimensional conformal manifold
$(\mathbf{M}^3,c)$ (also called the {\it twistor space} of {\bf M})
is a twistor space endowed with a contact structure, 
therefore we get (locally again) a self-dual manifold $(\mathbf{N}^4,c)$, in
which $(\mathbf{M}^3,c)$ is umbilic, and it is the {\it conformal infinity} of
an Einstein metric $g$ on $\mathbf{N}$ with non-zero scalar curvature
\cite{leb1}, \cite{leb2}.
\begin{defi}
Let $(\mathbf{M},c)$ be a conformal 3-manifold, that we shall suppose
{\em civilized} ({\em e.g.} 
geodesically connected for some metric in the conformal class%
).
The {\em LeBrun correspondence} associates to $\mathbf{M}$ the
(germ-unique) self-dual Einstein 4-manifold $\mathbf{N}$ such that the
twistor spaces of $\mathbf{M}$ and $\mathbf{N}$ coincide.
\end{defi}
\begin{prop}{\rm\cite{leb1}, \cite{leb2}}
In the LeBrun correspondence, $(\mathbf{M}^3,c)$ is an umbilic 
hypersurface of $(\mathbf{N}^4,c)$ (and has the induced conformal
structure) and the Einstein metric of $\mathbf{N}^4$ has a
second order pole at $\mathbf{M}^3$ (conformal infinity).
Conversely, in such a geometric setting, the twistor spaces of the
manifolds $\mathbf{M}$ and $\mathbf{N}$ coincide.
\end{prop}
\obs There is no {\it a priori} definition of a conformal infinity of
an open (real- or complex-) Riemannian manifold $X$, even if the metric is
complete (in the real framework). Here we consider uniquely
the case when this infinity is an (umbilic) submanifold (or boundary)
of a conformal extension of $X$, $\bar X\supset X$; the conformal
structure extends smoothly beyond the infinity. In other cases, in which
the notion of conformal infinity can still be defined, the conformal
structure is singular at infinity, which, in these cases, is no
longer conformal, but admits instead a $CR$ \cite{hitch},\cite{biq} or
a quaternionic contact structure \cite{biq}.

\section{Conformal infinity of a self-dual manifold}

The object of this section is to find a relationship between the
conformal invariants of a conformal infinity and of its self-dual
ambient arising from the LeBrun correspondence. The results are local,
and they hold in the complex as well as in the real Riemannian or in
the signature (2,2) pseudo-Riemannian framework. We begin by recalling
a few facts about the conformally invariant tensors in Riemannian
geometry. 
\medskip

For a $n$-dimensional Riemannian manifold $(M^n,g)$, the curvature has the 
following expression: 
\begin{eqnarray}\label{h3}
R^{ M}(X,Y)&=&(h\wedge\!\mbox{\bf\em I}\,)(X,Y)+W(X,Y),\mbox{ where}\\
(h\wedge\!\mbox{\bf\em I}\,)(X,Y) &:=& h(X,\cdot)\wedge Y-h(Y,\cdot)\wedge X,\;\forall X,Y\in
T M,
\end{eqnarray}
is {\it the suspension} by the identity \!{\bf\em I}\, of the {\it
  normalized Ricci tensor}
\begin{equation}\label{eq:rich}
  h=\frac{1}{2n(n-1)}\mbox{\it Scal}^g\cdot g+\frac{1}{n-2}\mbox{\it Ric}_0;
\end{equation}
 $Scal$ and $Ric_0$ are the scalar curvature, resp. the trace-free
 Ricci tensor, and, together with the Weyl tensor $W$, they are the
 irreducible components of the curvature under the orthogonal group if
 $n\ge 5$. If $n=3$, $W$ vanishes identically, and if $n=4$ it further
 decomposes in two irreducible components $W^+$, resp. $W^-$, called
 the {\it self-dual}, resp. {\it anti-self-dual} (or {\it
 positive}, resp. {\it negative})  Weyl tensor.

The Weyl tensor, viewed as a section in $\mbox{Hom}(\Lambda^2TM\otimes
TM,TM)$ (as a (3,1)--tensor), is conformally invariant, and, if $n\ge
4$, it completely determines, locally, the conformal structure of
$(M,[g])$ (for a self-dual manifold, $W^-\equiv 0$, thus the
Weyl tensor actually coincides with $W^+$). In dimension $n=3$, this function is
fulfilled by {\it the 
  Cotton-York tensor}, which can be defined in all dimensions by
\begin{equation}
  \label{eq:CY4}
  C(X,Y)(Z):=(\nb_Xh)(Y,Z)-(\nb_Yh)(X,Z),\;\forall X,Y,Z\in T M,
\end{equation}
and it can be shown that, for another metric $g':=e^{2\varphi} g$ in the
same conformal class, the corresponding Cotton-York tensor $C'$ is
related to $C$ by the formula \cite{gconf}:
\begin{equation}\label{cyweyl}
C'(X,Y)(Z)=C(X,Y)(Z)-\mbox{d}\varphi(W(X,Y)Z).
\end{equation}
In particular $C$ is conformally invariant along the zero set of $W$,
thus everywhere if $\dim M=3$.

\obs The Cotton-York tensor $C$ of $M$ is a 2-form with values in $T^*M$,
and it satisfies a {\it first} Bianchi identity,
as $h$ is a symmetric tensor, and also a {\it contracted}
(second) Bianchi identity, coming from the second Bianchi identity in
Riemannian geometry, \cite{gconf} :
\begin{eqnarray}\label{ar1}
\sum C(X,Y)(Z)&=&0 \;\mbox{ circular sum};\\
\label{ar2}
\sum C(X,e_i)(e_i)&=&0\;\mbox{ trace over an orthonormal basis}.
\end{eqnarray}
This means that $C$ is an irreducible tensor if $n=3$ or $n>4$,
and, if $n=4$, $C$ has two irreducible components, the {\it self-dual},
resp. the {\it anti-self-dual} Cotton-York tensor
$$C^+\in \Lambda^+M\otimes\Lambda^1 M,\ \mbox{resp.}\ C^-\in\Lambda^-
M\otimes\Lambda^1 M.$$
They both satisfy (\ref{ar1}) and (\ref{ar2}) (note that these two
relations are equivalent in their case). 


The Cotton-York tensor is related to the Weyl tensor of $M$ by the
formula \cite{gconf}:
\begin{equation}
  \label{eq:dw=cy}
  \delta W=C,
\end{equation}
where $\delta:\Lambda^2 M\otimes\Lambda^2 M\rightarrow
\Lambda^2 M\otimes\Lambda^1 M$ is induced by the codifferential on the
second factor, and by the Levi-Civita connection $\nb$. Then, again if
$n=4$, $C^+$
has to be the component of $\delta W$ in
$\Lambda^+ M\otimes\Lambda^1 M$, and we know that the restriction of
$W^-$ to $\Lambda^+ M\otimes\Lambda^2 M$ is identically zero. This means
that
\begin{eqnarray}
  \label{ar+}
  \delta W^+&=&C^+,\;\mbox{and also}\\
  \delta W^-&=&C^-.
\end{eqnarray}

We have thus:
\begin{lem}\label{C-}
On a self-dual manifold, $C^-$ vanishes identically. 
\end{lem}


We consider now the situation in the LeBrun correspondence~: Let
$(M,c)$ be a 3-dimensional conformal manifold, and we suppose, 
without any local loss of generality, that
it is the conformal infinity of the self-dual manifold $(N,c)$ (no use
will be made of the Einstein metric on it); $M\subset
N$ is, thus, an umbilic hypersurface, such that the restriction of the 
conformal structure $c$ of $ N$ to $ M$ is non-degenerate (equivalently,
$T M$ is nowhere tangent to an isotropic cone) and coincides with the
conformal structure, still denoted by $c$, on $ M$. 

If we introduce
the Hodge operator $*^{\!\scriptscriptstyle{M}}:\Lambda^2M\rightarrow
\Lambda^1 M$, then the curvature tensor $R^ M$ is equivalent to the
symmetric 2-tensor $*^{\!\scriptscriptstyle{ M}}\circ R^ M\circ
*^{\!\scriptscriptstyle{ M}}$. A straightforward application of the
above formula yields 
\begin{equation}\label{R*}
*R^ M:=*^{\!\scriptscriptstyle{ M}}\circ R^ M \circ
*^{\!\scriptscriptstyle{ M}}=-h+(\mathrm{tr}h)I. 
\end{equation}

For the 4-dimensional manifold $ N$, the components of the Riemannian 
curvature can also be expressed as eigenspaces of $*$-type operators.
Namely, considering $R:=R^ N$ as a symmetric endomorphism of $\Lambda^2 N
=\Lambda^+ N\oplus\Lambda^- N$, $W^+$ is the trace-free component of $R$
in $\mathrm{End}(\Lambda^+ N)$, and $W^-$ is the trace-free component
of $R$ in $\mathrm{End}(\Lambda^- N)$ \cite{st}. (The scalar curvature
is four times the trace of $R|_{\Lambda^+}$ or of $R|_{\Lambda^-}$,
and the trace-free Ricci tensor is identified to the component of $R$
sending $\Lambda^+$ into $\Lambda^-$ \cite{st}.)


We can canonically identify $\Lambda^+ N$ and $\Lambda^- N$, restricted
to $ M\subset  N$, to $\Lambda^2 M$, by:
\begin{equation}
  \label{eq:star}
  \begin{array}{lcr}
\Lambda^2 M\ni\alpha& \mapsto &\alpha+*^{\!\scriptscriptstyle{ N}}
\alpha\in\Lambda^+ N\\
\Lambda^2 M\ni\alpha& \mapsto &\alpha-*^{\!\scriptscriptstyle{ N}}
\alpha\in\Lambda^- N.
  \end{array}
\end{equation}

Our first result is:
\begin{thrm}\label{umbilic}
  Let $ M$ be an umbilic hypersurface of a self-dual manifold
  $ N$. Then:

(i) The Weyl tensor of $ N$ vanishes along $ M$~:
$$W^+|_M\equiv 0;$$

(ii) The Cotton-York tensor of $ M$ is related to the self-dual Weyl
tensor of $ N$ by the formula:
$$g(\nb_\nu W^+(A),B)_x=-C(A)(*^{\!\scriptscriptstyle{ M}}B)_x,\
\forall x\in M$$
where $A,B\in\Lambda^2 T_xM$, $\nu\perp T_x M$ is unitary for the metric
$g$, and the Hodge operator $ *^{\!\scriptscriptstyle{ M}}$ is induced
by $g$ and the orientation on $ M$ admitting $\nu$ as an exterior
normal vector.

(iii) The restriction to $M$ of the (self-dual) Cotton-York tensor of $N$ is
equal to the Cotton-York tensor of $M$~:
$$C^+(X,Y)(Z)=C^ M(X,Y)(Z),\;\forall X,Y,Z\in T M.$$
\end{thrm}
\begin{proof} The claimed identities are
  conformally invariant~: for {\it (i)} it is obvious, and the conformal invariance of {\it (iii)}
  follows from {\it (i)} and (\ref{cyweyl}); to see that for {\it
  (ii)}, let $X,Y,Z,\nu$ be a $g$-ortho\-normal
oriented basis of $ N$, such that $X,Y,Z$ is a $g$-orthonormal basis on $
M$ giving the orientation as above. Then $*^{\!\scriptscriptstyle{ M}}
(Z\wedge X)=Y$, and, if we take $A:=X\wedge Y,\, B:=Z\wedge X$, the
identity {\it (ii)} becomes
\begin{equation}
  \label{eq:CYW}
  \langle \nb_\nu W^+(X,Y)Z,X\rangle =-C(X,Y)(Y),
\end{equation}
where angle brackets denote the scalar product induced by $g$.
 
The tensors $W^+, C$, in the above form, are independent of the chosen
metric $g$ \cite{gconf}, which depends on the normal vector $\nu$,
supposed to be $g$-unitary. If $\nu':=\lambda\nu$, for
$\lambda\in\C^*$, then the corresponding metric $g'=\lambda^{-2}g$,
and also ${*^{\!\scriptscriptstyle{ M}}}'=\lambda^{-1}
*^{\!\scriptscriptstyle{ M}}$, thus the identity (\ref{eq:CYW}) for
$\nu',g'$ is equivalent to the one for $\nu,g$.

\obs As $W^+$ is the trace-free component of the Riemannian curvature
contained in $\mbox{\it End}(\Lambda^+ N)$, and is symmetric, it is enough
to evaluate it on pairs $A,B\in\Lambda^2 M\simeq\Lambda^+ N$ which are
unitary and orthogonal for the metric $g$, therefore the 
check of the equation (\ref{eq:CYW}) will prove the Theorem.
\medskip

As $W^\pm$ are $*^{\!\scriptscriptstyle{ N}}$-eigenvectors in
$\mbox{\it End}_0(\Lambda^2 N)$ (the space of trace-free endomorphisms of
$\Lambda^2 N$), they are determined by the following formulas, where
$X,Y,Z$ is any oriented orthonormal basis of $T M$:
\begin{eqnarray}
  \label{arw+}
 \  \langle W^+(X,Y)Z,X\rangle &\!\! =&\!\!\frac{1}{4}(\langle
 R(X,Y)Z,X\rangle +\langle R(Z,\nu)Y,\nu\rangle +\\ 
\nonumber  &&+\langle R(X,Y)Y,\nu\rangle + \langle R(Z,\nu)Z,X\rangle )\\
\label{arw-}
 \  \langle W^-(X,Y)Z,X\rangle &\!\! =&\!\!\frac{1}{4}(\langle
 R(X,Y)Z,X\rangle +\langle R(Z,\nu)Y,\nu\rangle -\\ 
\nonumber  &&-\langle R(X,Y)Y,\nu\rangle - \langle R(Z,\nu)Z,X\rangle ), 
\end{eqnarray}
where $X,Y,Z,\nu$ is supposed to be a local extension, around a region
of $ M$, of the $g$-orthonormal frame used in (\ref{eq:CYW}). As $ N$ is
self-dual, $W^-$ is identically zero, thus, in the
points $x\in  M$, we have  
\begin{equation}
  \label{eq:W1}
  \langle W^+(X,Y)Z,X\rangle _x=\frac{1}{2}(\langle R(X,Y)Y,\nu\rangle
  + \langle R(Z,X)Z,\nu\rangle )_x. 
\end{equation}

It is a standard fact that, if $ M$ is umbilic, there is a local metric
$g$ in the conformal class $c$ of $ N$, such that, for $g$, $ M$ is
totally geodesic. Without any loss of generality, because of the
conformal invariance of the claimed identities (see above), we fix
such a metric. Then we have 
\begin{equation}
  \label{eq:tgeod}
  R(X',Y')Z'=R^ M(X',Y')Z',\;\forall X',Y',Z'\in T M,
\end{equation}
which, together with (\ref{eq:W1}), implies that $W^+|_ M\equiv 0$, and
thus proves the point {\it (i)} in the Theorem.
\smallskip

On the other hand, (\ref{eq:tgeod}), together with (\ref{eq:W1}) and
(\ref{arw-}), yield
\begin{equation}
  \label{eq:rq}
  \langle R(X,Y)Z,X\rangle _x+\langle R(Z,\nu)Y,\nu\rangle _x=0,
  \forall x\in M.
\end{equation}

Let us compute now the normal derivative of $W^+$ in a point $x\in  M$;
we suppose that $X,Y,Z,\nu$ are locally extended by an orthonormal
frame, and that they are parallel at $x$ (we omit, for simplicity of
notation, the point $x$ in the following lines:
$$
\langle \nb_\nu W^+(X,Y)Z,X\rangle  =\frac{1}{2}(\langle \nb_\nu
R(X,Y)Y,\nu\rangle +\langle \nb_\nu 
R(Z,X)Z,\nu\rangle ),$$
from (\ref{eq:W1}). This is then equal to:
$$\begin{array}{rcl}
\langle \nb_\nu W^+(X,Y)Z,X\rangle &\!\!\! =&\!\!\!-\frac{1}{2}(\langle \nb_X
R(Y,\nu)Y,\nu\rangle +\langle \nb_Y
R(\nu,X)Y,\nu\rangle +\\
 &&+\langle \nb_Z R(X,\nu)Z,\nu\rangle  +\langle \nb_X R(\nu,Z)Z,\nu
 \rangle ),\end{array}$$  
from the second Bianchi identity. Then we have
$$\begin{array}{rcl}
\langle \nb_\nu W^+(X,Y)Z,X\rangle &\!\!\! =&
\!\!\!\frac{1}{2}(\langle \nb_X R(Z,X)Z,X\rangle +\langle \nb_Y R(Z,Y)Z,X\rangle +\\
&&+\langle \nb_Z R(Y,Z)X,Y\rangle +\langle \nb_X R(Y,X)X,Y\rangle )\end{array}$$
from analogs of (\ref{eq:rq}). Then
$$\begin{array}{rcl}
\langle \nb_\nu W^+(X,Y)Z,X\rangle &\!\!\!\!\! =&\!\!\!\!\frac{1}{2}(\langle \nb_X
R^ M(Z,X)Z,X\rangle  +\langle \nb_Y R^ M(Z,Y)Z,X\rangle  +\\
&&\!\!\! +\langle \nb_Z R^ M(Y,Z)X,Y\rangle  +\langle \nb_X
R^ M(Y,X)X,Y\rangle ) 
\end{array}$$
from (\ref{eq:tgeod})
$$\begin{array}{rcl}
\langle \nb_\nu W^+(X,Y)Z,X\rangle &\!\!\! =&\!\!\!\frac{1}{2}(\nb_X
h(Z,Z)+\nb_X h(X,X)+\nb_Y h(Y,X)-\\ 
&&-\nb_Z h(Z,X)-\nb_X h(X,X)-\nb_X h(Y,Y)),
\end{array}$$
from (\ref{h3}). Finally, from (\ref{ar2}), we get
$$
\langle \nb_\nu W^+(X,Y)Z,X\rangle =\frac{1}{2}(C(X,Z)(Z)-C(X,Y)(Y))=-C(X,Y)(Y)$$
This proves equation (\ref{eq:CYW}) and the point {\it (ii)} in the
Theorem.

For the point {\it (iii)}, we use (\ref{eq:CYW}) and (\ref{ar+});
from the codifferential of $W^+$, only the derivative along the normal
vector, $\nu$, can be non-zero, as $W^+$ vanishes along $M$. This
proves the Theorem.
\end{proof}

\obs The point {\it (i)} gives a condition for a self-dual manifold to
admit an umbilic hypersurface~: $W^+$ has to vanish along it,
generically at order 0 (following {\it (ii)}), thus such a
hypersurface, if it exists, is locally defined as the zero set of
$W^+$.

Considering the umbilic hypersurfaces which arise from the LeBrun
correspondence, the point {\it (i)} gives a condition for an open
self-dual manifold to admit such a conformal infinity, namely it has to be
asymptotically conformally flat.

\section{Null-geodesics of complex conformal manifolds}
\subsection{Properties of the twistor space of a 3-dimensional
  conformal manifold. }
Consider ({\bf M},$c$) a complex 3-dimensional conformal manifold. In
some topological conditions ({\bf M} has to be {\it civilized}
\cite{leb1}~; as a geodesically convex set is always of this type
\cite{leb5}, any point has a basis of civilized neighbourhoods), the
{\it twistor space} of {\bf M} is defined as the 
space $Z$ of null-geodesics of {\bf M}, and it is a complex
3-manifold, containing {\it twistor lines} ({\it i.e.} rational curves
with normal bundle isomorphic to $\o1\oplus\o1$), and endowed with a
distribution of 2-planes $F_{\bar\gamma}\subset T_{\bar\gamma}Z$ which
is a contact structure \cite{leb1}. 

We denote by $\bar\gamma$ the point of $Z$ corresponding, in the
following way, to the
null-geodesic $\gamma$ :
some twistor lines tangent to $F_{\bar\gamma}$ (actually a non-empty open set
in the space of curves of $Z$, tangent to $F$) 
correspond to the points
of the null-geodesic $\gamma$\,\footnote{deformations of these twistor
  lines are, in general, not tangent to the distribution $F$; they
  correspond to points in the self-dual ambient $N$ (and outside of
  $M$) arising from the
  LeBrun correspondence.}.

The first question we raise is whether there exist twistor lines
tangent to any direction of a given 2-plane $F_{\bar\gamma}$; the
answer is:

\begin{thrm}\label{compact3}
  Let $Z$ be a twistor space of a conformal civilized 3-manifold {\bf M}; let
  $F_{\bar\gamma}\subset T_{\bar\gamma}Z$ be its contact
  structure. Suppose there is a point $\bar\gamma\in Z$ such that
  there are twistor lines tangent
  to any direction in $F_{\bar\gamma}$. Then $Z$ is
  projectively flat, and {\bf M} is conformally flat.
\end{thrm}
This follows directly from \cite{FW}, Theorem 3, which has a similar
statement referring to the twistor space of a self-dual manifold; we
simply apply it to the ambient {\bf N} from the LeBrun correspondence; its
conformal flatness implies the flatness of {\bf M} (Theorem
\ref{umbilic}).

\obs The key point in the above cited Theorem is a {\it twistorial
  interpretation} of the Weyl tensor of a self-dual manifold {\bf N}
  \cite{FW}, Theorem 2, together with the remark that, for a given
  2-plane $F$ in $T_{\bar\beta}Z$, the union of all twistor lines tangent
  to it (supposing there exists one pointing in any direction of $F$)
  is a complex surface which is {\it smooth} at $\bar\beta\in Z$. The
  above cited Theorem and the following one (Theorem $3'$ from
  \cite{FW}) show that this situation implies the vanishing of the
  Weyl tensor of {\bf N}, $W^+$, in certain directions:
  
\begin{thrm}\label{3'}{\rm \cite{FW}}
Let $Z$ be the twistor space of the (civilized) self-dual manifold {\bf
  N}, and let $\bar\beta\in Z$ be a point in $Z$, corresponding to the
$\beta$-surface $\beta\subset${\bf N}; let $F^\gamma\subset T_{\bar\beta}Z$
be a 2-plane, corresponding to the null-geodesic
$\gamma\subset\beta$ \cite{leb3}. Suppose that, for each direction
$\sigma\in\P(F^\gamma)\subset\P(T_{\bar\beta} Z)$, there is a smooth
(non-necessarily compact) curve $Z_\sigma$ tangent to $\sigma$, such
that: 

(i) if $\sigma$ is tangent to a twistor line $Z_x$ at $\bar\beta$, then
$Z_\sigma=Z_x$;

(ii) $Z_\sigma$ varies smoothly with $\sigma\in\P(F^\gamma)$.

Then 
$$\bar Z_\beta^\gamma:=\bigcup_{\sigma\in\P(F^\gamma)}Z_\sigma$$ is a smooth
surface around $\beta$, and
$W^+(F^\gamma_x)=0,\,\forall x\in\gamma$, where $F^\gamma_x\subset
T_x \mathbf{M}$ is the $\alpha$-plane containing $\dot\gamma$.
\end{thrm}
 In other words, if the {\it integral $\alpha$-cone} corresponding to
 the 2-plane $F\subset T_{\bar\gamma}Z$ --- defined as the
 union of all twistor lines tangent $F$ \cite{FW} --- can be
 completed to a surface, smooth around its ``vertex'' $\bar\gamma$, then $W^+$
 vanishes along the null-geodesic $\gamma$ (whose points correspond to
 the twistor lines that constitute the $\alpha$-cone \cite{FW}).

\subsection{Compact null-geodesics and conformal flatness. } Our main
 result is:
\begin{thrm}
  \label{comp3} Let {\bf M} be a conformal $n$-manifold containing an
  immersed rational curve as null-geodesic. Then {\bf M} is conformally flat.
\end{thrm}
This fact has been proven by Ye \cite{ye} for complex projective
manifolds --- we are grateful to B. Klingler for having brought this
paper into our attention. In the above Theorem, we do not make any assumption on the
topology of the manifold, but only of one null-geodesic contained in it.

\begin{proof} The proof is different in the cases $n>3$ and $n=3$; one
  of the reasons is that conformal flatness reduces, in higher
  dimensions, to the vanishing of the Weyl tensor, while in dimension
  3 it is a higher-order condition more difficult to handle. The first
  step (common to all cases) is to prove that a small deformation (seen just as a 
compact submanifold of {\bf M}, \cite{kod}) of such a compact
null-geodesic $\gamma$ is still a compact null-geodesic, and to characterize
  the global sections of the normal bundle of $\gamma$ as locally determined
  by {\it Jacobi fields}. 

\begin{lem}\label{lema1} Let $\gamma$ be a (immersed) null-geodesic in
  $(\mathbf{M},c)$. Let $J$ be a
  vector field along $\gamma$. Then the condition $\dot
  J\perp\dot\gamma$ (where $\dot J:=\nabla^g_{\dot\gamma}J$) is
  independent of the metric $g$ with respect to which we take the
  derivative $\nabla^g$.
\end{lem}
\begin{proof} The relation between two Levi-Civita connections (or,
  more generally, {\it Weyl structures}) of metrics in the same
  conformal class, is given by \cite{gconf}:
\begin{equation}\label{thetaw}
\nabla'_XY-\nabla_XY=\theta(X)Y+\theta(Y)X-\theta^\sharp g(X,Y),
\end{equation}
where $\theta$ is a 1-form, and the rising of indices in
$\theta^\sharp$ is made using the same (arbitrary) metric $g\in c$ as
in the scalar product $g(X,Y)$. The Lemma immediately follows.
\end{proof}

Denote by $N(\gamma)$ the normal bundle of $\gamma$ in {\bf M}, and by
  $N^\perp(\gamma)$ its subbundle 
  represented by vectors orthogonal to $\dot\gamma$. Fix a metric $g$
  in the conformal class $c$. Let $J$ be a
  Jacobi field along $\gamma$, satisfying to the Jacobi equation
$$\ddot J=R^g(\dot\gamma,J)\dot\gamma.$$
It represents an infinitesimal deformation of $\gamma$ through
null-geodesics if and only if $\dot J\perp\dot\gamma$. $J$ induces a
section in the normal bundle $N(\gamma)$, or $N^\perp(\gamma)$ if
$J\perp\dot\gamma$ in a point, hence everywhere. We want to prove that
this section is independent of the connection $\nabla^g$:
\begin{prop}\label{prop1}
  The Jacobi equations on a null-geodesic $\gamma\subset${\bf M}
  induce a second order linear 
  differential operator $P$ on $N(\gamma)$, which, restricted to the
  sections $J$ such that $\dot J\perp\dot\gamma$, depends only on the
  conformal structure $c$ of {\bf M}. In particular, $P$ restricted to
  $N^\perp(\gamma)$ is conformally invariant.
\end{prop}
\begin{proof} For a Levi-Civita connection $\nb$ of a local metric on
  {\bf M}, we locally define the following differential operator on
  the sections of $T\mathbf{M}|_\gamma$:
$$P:\Gamma(T\mathbf{M}|_\gamma\otimes
S^2(T\gamma))\rightarrow\Gamma(T\mathbf{M}|_\gamma),$$
by $P(Y;X,X):=\nb_X\nb_XY-\nb_{\nb_XX}Y-R(X,Y)X$. Because $\gamma$ is
a null-geodesic, $P$ induces
a (local) differential operator on $N(\gamma)$, and we need to
relate $P$ to the corresponding
operator $P'$ induced by another connection
$\nabla'$. First we write 
$$P(Y,X,X)=\nb_X[X,Y]+\nb_{[X,Y]}X-[\nb_XX,Y],$$ 
then we recall that another Levi-Civita connection $\nb'$ is related
to $\nb$ by the formula (\ref{thetaw}), such that we get
$$P'(Y;X,X)-P(Y;X,X)=2[Y.(\theta(X))-\theta([X,Y])] X
-g(\nabla_XY,X)\theta^\sharp,$$
which is identically zero modulo $T\gamma$, provided that
$\nabla_XY\perp X$ (the latter condition being independent of the
Levi-Civita connection, according to the previous Lemma).
\end{proof}
Using the fact that $\cp1$ is the union of two contractible sets
$U_1\cup U_2$ (on each of
which the Jacobi equation, with any initial condition --- the same for
$U_1$ and for $U_2$ ---
in $x_0\in U_1\cap U_2$, has a unique solution --- and these solutions
necessarily coincide on the connected intersection $U_1\cap U_2$), 
we immediately get:
   \begin{prop}\label{p5}
     Let $\gamma$ be an immersed null-geodesic, diffeomorphic to a projective line
     $\cp1$. Then any local Jacobi field $J$ with $\dot
     J\perp\dot\gamma$ induces a global normal field $\nu^J$ on
     $\gamma$.
   \end{prop}
This has important consequences about the normal bundle of $\gamma$ in
{\bf M}, as Jacobi fields provide it with global sections; in particular,
$N(\gamma)/N^\perp(\gamma)$ is a line bundle admitting
nowhere-vanishing sections, hence it is trivial; on the other hand,
$N^\perp(\gamma)$ is a $(n-2)$--rank bundle over $\cp1$, admitting
sections with any prescribed 1-jet (induced, again, by some appropriate
Jacobi fields), hence 
\begin{equation}\label{fibnorm}
N^\perp(\gamma)\simeq \bigoplus_{k=1}^{n-2}\o{a_k};\ N(\gamma)\simeq
\bigoplus_{k=1}^{n-2}\o{a_k}\oplus\o0;\ a_k\in\N^*.
\end{equation}
For $a_k\in\Z$, this is the general form of a vector bundle over
$\cp1$, according to a theorem of Grothendieck; the condition $a_k\ge
1$ arises from the existence of sections of $N^\perp(\gamma)$ with
prescribed 1-jet.
We are going to show later that all $a_k$ are equal to 1.
First we prove:
\begin{prop}\label{tub3}
Null-geodesics close to a compact, simply-connected one are also
compact and simply-connected, and they are generically embedded.
\end{prop}
\begin{proof} We consider the projectivized bundle $\P(\mathcal{C})$ of
  the isotropic cone $\mathcal{C}\subset T\mathbf{M}$. It is a standard fact
  \cite{leb1} that any null-geodesic $\gamma\subset \mathbf{M}$ has a canonical
  {\it horizontal} lift $\tilde\gamma\subset\P(\mathcal{C})$ (depending only
  on the conformal structure), such that
  $\pi_*(T_s\tilde\gamma)=T_x\gamma$, where
  $\pi:\P(\mathcal{C})\longrightarrow \mathbf{M}$ is the projection, and
  $s\in\pi^{-1}(x)$.

 Note that $\tilde\gamma$ is always embedded, even if $\gamma$ may have
 self-intersections (it is always immersed). 

The lifts of the null-geodesics of {\bf M} consist in a foliation of
$\P(\mathcal{C})$, which has a compact, simply-connected leaf, namely
the lift $\tilde\gamma$ of our compact, simply-connected null-geodesic
$\gamma$. By Reeb's stability Theorem \cite{reeb}, then there is a
{\it saturated} tubular neighbourhood of $\tilde\gamma$, diffeomorphic to
$\tilde\gamma\times D$ (where $D\subset\C^{n-2}$ is a polydisc), such
that the leaves close to $\tilde\gamma$ are 
identified, via the above diffeomorphism, to the (compact and
simply-connected) curves $\tilde\gamma\times\{\mbox{z}\},\ z\in D$.

So all null-geodesics close to $\gamma$ are compact and simply
connected. If $\gamma$ has self-intersections at the points
$x_1,\dots,x_k$, we blow-up {\bf M} at those points, and the lift of
$\gamma$ is now embedded. So must be then the lifts of the
null-geodesics close to $\gamma$, as they are now deformations of the
lifted (hence, embedded) curve. But, generically, such curves avoid
the finite set of points $x_1,\dots,x_k$; the corresponding
null-geodesics must have been embedded from the beginning.
\end{proof}

From now on, according to the previous Proposition, we may suppose
that $\gamma$ is a compact, simply-connected, {\it embedded}
null-geodesic.
\medskip

We compute the normal bundle of $\gamma$ in {\bf M}, using the relation
(\ref{fibnorm}) and the projection $\pi:\P(\mathcal{C})\longrightarrow
\mathbf{M}$, as follows:
We have the following exact sequence of bundles:
\begin{equation}
  \label{eq:opi}
  0\rightarrow N^\pi(\tilde\gamma)\rightarrow
  N(\tilde\gamma)\rightarrow\pi^* N(\gamma)\rightarrow 0,
\end{equation}
where $N^\pi(\tilde\gamma)$ is the normal subbundle of $\tilde\gamma$
represented by vectors tangent to the fibers of $\pi$ and
$N(\tilde\gamma)$ is the normal bundle of $\tilde\gamma$ in
$\P(\mathcal{C})$. In a point
$T_x\gamma\in\tilde\gamma\subset\P(\mathcal{C})$, the fiber of $\pi$ is equal to
$\P(\mathcal{C})_x$, so the tangent space to it is isomorphic to
$\mbox{Hom}(T_x\gamma,N_x^\perp(\gamma))$, for the projective
variety $\P(\mathcal{C})_x\subset\P(T_x\mathbf{M})$. Thus 
$$N^\pi(\tilde\gamma)\simeq\mbox{Hom}(T\gamma,N^\perp(\gamma))\simeq
\o {-2}\otimes N^\perp(\gamma),$$
as $T\gamma\simeq T\cp1\simeq\o{-2}$.

The central bundle in the exact sequence (\ref{eq:opi}) is
trivial, because $\tilde\gamma$ is a leaf of a foliation. 
(\ref{fibnorm}) and (\ref{eq:opi}) imply that the Chern numbers
$a_1,\dots,a_{n-2}$ are subject to the following constraint: 
$$\sum_{k=1}^{n-2}(2a_k-2)=0,$$
thus, as $a_k\ge 1$, we have $a_k=1,\;\forall k=\overline{1,n-2}$.
We have then:
\begin{prop}\label{deform}
The normal bundle of a compact, simply-connected, null-geodesic
  $\gamma$ in {\bf M} is isomorphic to 
$$N(\gamma)\simeq\left(\C^{n-2}\otimes\o1\right)\oplus\o0,$$
and all its global sections are induced by Jacobi fields $J$ such that
$\dot J\perp\dot\gamma$. Moreover, the deformations of $\gamma$ as a
compact curve coincide with the null-geodesics close to $\gamma$.
\end{prop}
The last statement follows from the expression of the normal bundle,
and a Theorem of Kodaira \cite{kod}: the normal bundle satisfies
$H^1(N(\gamma))=0$, thus the space of deformations of $\gamma$ as a
compact curve has dimension equal to $\dim
H^0(N(\gamma))=\dim\Gamma(N(\gamma))=2n-3$, which 
is precisely the dimension of the space of null-geodesics, defined
locally, over a geodesically convex open set, as the space of the
leaves of the horizontal foliation of $\P(\mathcal{C})$ \cite{leb5}.
We conclude using the fact that all null-geodesics close to $\gamma$
are deformations of this one (as a compact, and simply-connected,
curve).
\medskip

We return to the proof of Theorem \ref{comp3}. Consider first the case
when the dimension of {\bf M}, $n>3$. We are going to show that the
Weyl tensor of {\bf M}, $W$, is identically zero (a special 
sub-case is $n=4$, when $W=W^+ + W^-$). For simplicity, suppose first
that $n>4$, and consider the fiber of $N^\perp(\gamma)$ at an
arbitrary point $x\in\gamma$: it has a
non-degenerate conformal structure, induced from {\bf M}, and the
isotropy cone spans the whole fiber $N^\perp(\gamma)_x$ (this still
holds for $n=4$, but not for $n=3$). Let $L\subset N^\perp(\gamma)$ be
an isotropic line. We have:
\begin{lem}\label{paralel} Let $\gamma^U$ be an open set of the
  null-geodesic $\gamma$, on which local metrics $g,g'\in c$ are well defined. 
If a (locally defined, over $\gamma^U$) line subbundle $L\subset
N^\perp(\gamma)$ is parallel (or stable) for 
$\nabla^g$, then it is parallel for $\nabla^{g'}$ as well.
\end{lem}
The proof is a straightforward application of (\ref{thetaw}).
\smallskip

Let $(L_1)_x,\dots,(L_{n-2})_x$ be linearly independent isotropic lines
in $N^\perp(\gamma)_x$. According to the previous Lemma, and to the
fact that $\gamma$ is simply-connected, their parallel transport over
$\gamma$ does not depend on any Levi-Civita connection of a metric in
the conformal class. We get, thus, a {\it global} splitting 
\begin{equation}\label{split}
N^\perp(\gamma)=L_1\oplus\dots\oplus L_{n-2},
\end{equation}
where the line bundles $L_i,\; i=\overline{1,n-2}$ are all isotropic
and parallel.

All these bundles are isomorphic to $\o {b_i},\; b_i\in\Z$. As their
sections are also sections of $N^\perp(\gamma)\simeq
\C^{n-2}\otimes\o1$, they cannot vanish at more that 1 point, for each
$L_i$, thus $b_i\le 1,\; i=\overline{1,n-2}$. On the other hand, the sum of all
$b_i$'s has to be $n-2$, thus $b_i=1,\; \forall i=\overline{1,n-2}$.

Let $\phi_i$ be a section of $L_i$; from Proposition \ref{deform},
it is locally represented by a 
Jacobi field $J_i$, for the metric $g\in c$. From the Jacobi
equation, by taking the scalar product with $J$, we get:
\begin{equation}
  \label{eq:curv}
  g(R(\dot\gamma,J)\dot\gamma,J)=0,
\end{equation}
  and it is easy to see that, because of the fact that all scalar
  products involving $\dot\gamma$ and $J$ are 0, the term
  $h\wedge\!\mbox{\bf\em I}\,$ of the curvature (\ref{h3}) satisfies
  the above relation identically. 
 This equation holds, at $x$, for any isotropic vector
 $J_x\perp\dot\gamma_x$, but we may consider also other compact,
 simply-connected, null-geodesics $\gamma'$, containing $x$, and close
 to $\gamma$ (namely, small deformations of the compact curve $\gamma$).

For any 2-plane $F\subset T_x${\bf M}, we denote by
$R^F$ the {\it sectional curvature} of $F$:
 $$R^F:S^2(\Lambda^2F)\rightarrow\C,\  R^F(X\wedge Y,X\wedge
 Y):=\langle R(X\wedge Y),X\wedge Y\rangle,\; \forall X\wedge Y\in\Lambda^2F,$$
 and we have seen that, if $F$ is totally isotropic, $R^F$ depends
 only on $W$ (and on the metric $g$ only via the scalar product
 $\langle.,.\rangle$). 

\begin{lem}\label{weylF}
If $\dim\mathbf{M}>4$, the Weyl tensor at $x\in\mathbf{M}$, $W_x$, is
determined by the sectional curvatures $\{R^F,\ F\in
U(F_0)\}$, for $U(F_0)$  a small neighbourhood of the totally isotropic
arbitrary 2-plane $F_0\subset T_x\mathbf{M}$ in the Grassmanian of
totally isotropic 2-planes at $x$. 
\end{lem}
\obs A similar statement holds in dimension 4, but in that
case, the Grassmanian of totally isotropic 2-planes has 2 connected
components; as a consequence, $W^+$ is determined by the sectional
curvatures of $\alpha$-planes, and $W^-$ by the sectional curvatures
of the $\beta$-planes \cite{FW}, Proposition 2.
\begin{proof} This reduces to the claim that $W_x=0$ if and only if 
\begin{equation}\label{www}
R^F=0,\  \forall F\in U(F_0),
\end{equation}
which is a problem of linear algebra. If we
  consider the space $\mathcal{K}$ of all curvature tensors $R'$ satisfying
  $(R')^F=0,\ \forall F\in U(F_0)$, then this is a vector space, which
  is invariant to the action of $\mathfrak{so}(n,\C)$ (which is the
  Lie algebra infinitesimal action corresponding to the action of
  $SO(n,\C)$ --- note that the Grassmanian of totally isotropic planes
  is preserved by this action). But there are only 3
  $\mathfrak{so}(n,\C)$-irreducible components of the space of
  curvature tensors, and we have seen that for the Ricci-like tensors
  $h\wedge\!\mbox{\bf\em I}\,$, the totally isotropic planes always
  have zero sectional curvature. Then either any Weyl tensor
  satisfying (\ref{www}) is zero at $x$, or $\mathcal{K}$ contains
  the whole space of curvature tensors. The latter possibility can
  easily be excluded by an example of a curvature tensor $K$ satisfying:
$$K(X_0,Y_0)X_0=A_0,\mbox{ where } \langle Y_0,A_0\rangle=1, \mbox{
  and }$$
$$ \langle X_0,X_0\rangle=\langle
X_0,Y_0\rangle=\langle Y_0,Y_0\rangle=\langle X_0,A_0\rangle=0.$$
\end{proof}
From (\ref{eq:curv}) and the previous Lemma we conclude that $W_x=0$
for any $x$ contained in a compact, simply-connected, null-geodesic;
but we know from Proposition \ref{deform} that the set of such points
contains  a neighbourhood of 
$\gamma$, thus, by analyticity of $W$, it vanishes identically.
\medskip

The proof is similar in dimension 4 (note that, in the self-dual case,
we can retrieve the result by applying Theorem \ref{3'}; this is how
we shall proceed for the case of dimension 3, using the LeBrun
correspondence); the difference with the higher-dimensional case is
that the splitting (\ref{split}) is canonical, $L_1$ corresponding,
say, to the $\alpha$-plane $F^\alpha$ containing $\dot\gamma$, and $L_2$ to the
$\beta$-plane $F^\beta$ containing $\dot\gamma$. It is important now that {\it
  each of} $L_1,L_2$ is isomorphic to $\o1$, because the vanishing of
$R^{F^\alpha}$ implies $W^+\equiv 0$, and the vanishing of
$R^{F^\beta}$ implies $W^-\equiv 0$ \cite{FW}, Proposition 2. Thus the
manifold $(\mathbf{M}^4,c)$ is conformally flat.
\medskip

Consider now the particular case where $n=3$. We are going to use the
LeBrun correspondence, then Theorem \ref{3'}, to prove that {\bf M} is then
conformally flat. Note that we cannot use directly Theorem \ref{umbilic} and the
  above proven result for self-dual manifolds, as the ambient
self-dual manifold $ N$ can only be defined for a {\it civilized}
(e.g. geodesically convex) 3-manifold.
\smallskip

We cover $\gamma$ with geodesically convex open sets $U_i, \,
i=\overline{1,n}$, such that:
\begin{equation}
  \label{eq:gcx3}
  \forall i\neq j \mbox{ such that }U_i\cap U_j\cap\gamma\neq\emptyset,
  \ \exists U_{ij}\supset(U_i\cup U_j),
\end{equation}
where $U_{ij}$ is still geodesically convex (with respect to some
particular Levi-Civita connection). This is possible by choosing $U_i,
\, i=\overline{1,n}$, small enough \cite{ehr}. Then we choose a relatively
compact tubular neighbourhood $N(r_0)$ of $\gamma$, such that its
closure is covered by the $U_i$'s. We may choose this tubular
neighbourhood small enough to be contained in the projection $U$, from
$\P(\mathcal{C})$, of a saturated neighbourhood (see Proposition
\ref{tub3}) of the lift $\tilde\gamma$.

We consider then the twistor spaces $Z_i$, the spaces of
null-geodesics of $U_i$. The compact, simply-connected, null-geodesics
close to $\gamma$ identify (diffeomorphically) the neighbourhoods of
$\bar\gamma^i\in Z_i$ with the space $Z$ of the deformations
(contained in $U$) of $\gamma$ as a compact curve. We can see then (a
small open set of) $Z$ as an open set common to all the $Z_i$'s:
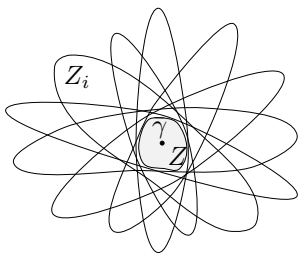
\begin{center}\input{weyl/3tw.pstex_t}\end{center}
Following LeBrun, we define the self-dual manifolds $ \mathbf{N}_i$ as
the spaces 
of projective lines in $Z_i$, with normal bundle $\o1\oplus\o1$. Then
$U_i$ is an umbilic hypersurface in $ \mathbf{N}_i$. 

The local twistor spaces $Z_i$ admit contact structures, which
coincide on $Z$, and contain projective lines $Z^i_x$ corresponding to
points $x\in\gamma\cap U_i$. If we denote by $Z_{ij}$ the twistor
space of $U_{ij}$, then $Z_{ij}$ is identified to an open set in $Z_i$
and, at the same time, to an open set in $Z_j$, in particular the
twistor lines $Z^i_x\subset Z^i$ and $Z^j_x\subset Z^j$ (corresponding
to the same point $x\in U_{ij}$) are identified. Thus $Z^i_x\cap Z$
and $Z^j_x\cap Z$ coincide and we denote by $Z_x$ this
(non-compact) curve in $Z$, and by $F$ the canonical contact structure
of $Z$ (restricted from the ones of $Z_i$).

\obs We already have obtained that the integral $\alpha$-cone ({\it
  i.e.} the union of twistor lines passing through $\bar\gamma$ and
  tangent to $F_{\bar\gamma}$, see the comment after Theorem \ref{3'})
  corresponding to $F_{\bar\gamma}$ is a part of a smooth surface: 
the union of the lines
$Z_x$, $x\in\gamma$. Thus, from Theorem \ref{3'}, the Weyl tensor $W_i^+$ of
the self-dual manifold $ \mathbf{N}_i$ vanishes on the $\alpha$-planes generated by
$T\gamma$. But this is nothing new: we know, from Theorem
\ref{umbilic}, that $W_i^+$ vanishes on $U_i$.
\smallskip

We intend to apply Theorem \ref{3'} to prove that $W_i^+$ vanishes on
points close to $U_i$, but generically in $ \mathbf{N}_i\smallsetminus
U_i$. We do that by 
showing that the integral $\alpha$-cones corresponding to planes $F'\subset
T_{\bar\gamma}Z$ close to $F$ are parts of smooth surfaces, then we
conclude using Theorem \ref{3'}.
\smallskip

First we choose Hermitian metrics $h_i$ on $Z_i$, such that they
coincide (with $h$) on $Z$. We have a diffeomorphism between $\gamma$
and $\P(F_{\bar\gamma})$, so we choose relatively compact open sets in
$\P(T_{\bar\gamma}Z)$, covering $\P(F_{\bar\gamma})$, with the
following properties: As the metrics $h_i$ induce metrics on 
$\mathbf{N}_i$, we first choose a small enough distance $r_1>0$ such
that
\begin{enumerate}
\item $\forall i$, there is a sub-covering $V_i\Subset U_i$ of
$\gamma$ such that the ``tubular neighbourhoods'' $Q_i:=\{y\in
\mathbf{N}_i\,|\,\mathrm{d}(y,\bar V_i)\leq r_1\,\, \pi_i(y)\in \bar
V_i\cap\gamma\}$ are compact ($\mathrm{d}(y,\bar V_i)$ is the distance
between $y$ and $\bar V_i$, and $\pi_i$ is the ``orthogonal projection''
--- for the Hermitian metric --- from $ \mathbf{N}_i$ to $\gamma\cap
U_i$; it is well defined because of the condition below);
\item $r_1$ is less than the bijectivity radius of the (Hermitian)
exponentials for the points of $\bar V_i$ in $ \mathbf{N}_i$, and for
the points of $\overline{V_i\cup V_j}$ in $ \mathbf{N}_{ij}$ (if
$U_i\cap U_j\cap\gamma\neq\emptyset.$).
\end{enumerate}
We have then
\begin{lem} For any $y_i\in Q_i\subset  \mathbf{N}_i$, $y_j\in
  Q_j\subset  \mathbf{N}_j$ 
  such that the curves $Z_{y_i}:=Z^i_{y_i}\cap Z$,
  $Z_{y_j}:=Z^j_{y_j}\cap Z$ are tangent to the same direction in
  $\bar\gamma\in Z$, the respective curves $Z_{y_i}, Z_{y_j}$ coincide.
\end{lem}
\begin{proof}
We first note that the projection $\pi_i$ from $ \mathbf{N}_i$ to
$\gamma\cap U_i$ is
equivalent to the $h$--orthogonal projection of the direction of
$T_{\bar\gamma}Z_{y_i}$ to a direction in $F_{\bar\gamma}$, so
$\pi_i(y_i)=\pi_j(y_j)=:y\in\gamma$; thus $y$ belongs to both $U_i$
and $U_j$, and we use again the twistor space $Z_{ij}$ to conclude
that $Z_{y_i}$ and $Z_{y_j}$ are ``restrictions'' to $Z$ of the same
projective line (as they both have the same tangent space at
$\bar\gamma$) $Z^{ij}_{y_{ij}}$, for a point $y_{ij}\in
\mathbf{N}_{ij}$.
\end{proof} 
 
Now we have a tubular neighbourhood $S\subset \P(T_{\bar\gamma}Z)$ of
$\P(F_{\bar\gamma})$, of radius $r_1/2$, such that, for any 2-plane
$F'\subset S$, the conditions in Theorem \ref{3'} are satisfied
(considering any of the local twistor spaces $Z_i$).

We recall that, via the LeBrun correspondence, a point $\bar\gamma_0$
in the twistor space of $\mathbf{M}_0$, $Z_0$, is identified to the
point $\bar\beta_0$ in the twistor space of $\mathbf{N}_0$, still
denoted by $Z_0$. They correspond to the null-geodesic
$\gamma_0\subset\mathbf{M}_0$, resp. to the $\beta$-surface
$\beta_0\subset\mathbf{N}_0$, such that $\gamma_0\subset\beta_0$
\cite{leb1}, \cite{FW}. The planes $F'$ above are included in
$T_{\bar\gamma}Z\equiv T_{\bar{\beta^i}} Z_i$, and they correspond to
null-geodesics in $\mathbf{N}_i$ contained in $\beta^i$ \cite{leb5},
\cite{FW}.

By Theorem \ref{3'}, we conclude that the Weyl tensor $W_i^+$ of $
\mathbf{N}_i$ vanishes along all null-geodesics of $ \mathbf{N}_i$,
close (in $Q_i$) to
$\gamma$ and included in the $\beta$-surface $\beta^i$, determined by
$\gamma$. This means that $W^+$ vanishes everywhere on $\beta^i$. By
deforming $\gamma$, we obtain that $W^+_i$ vanishes on a neighbourhood
of $U_i$ in $ \mathbf{N}_i$, hence $ \mathbf{N}_i$ is conformally flat.

It follows from Theorem \ref{umbilic} that $U_i$, hence {\bf M}, is
conformally flat.
\end{proof}

\bigskip
\begin{center}
{\sc {Mathematisches Institut\\ Humboldt Universit\"at zu Berlin\\
Unter den Linden 6, 10099 Berlin\\Germany}}\\e-mail: {\tt
belgun\@@mathematik.hu-berlin.de}
\end{center}
\bigskip
\end{document}

%% file: weyl/3tw.pstex_t
\begin{picture}(0,0)%
\epsfig{file=weyl/3tw.pstex}%
\end{picture}%
\setlength{\unitlength}{0.00022700in}%
\begingroup\makeatletter\ifx\SetFigFont\undefined
\def\x#1#2#3#4#5#6#7\relax{\def\x{#1#2#3#4#5#6}}%
\expandafter\x\fmtname xxxxxx\relax \def\y{splain}%
\ifx\x\y   
\gdef\SetFigFont#1#2#3{%
  \ifnum #1<17\tiny\else \ifnum #1<20\small\else
  \ifnum #1<24\normalsize\else \ifnum #1<29\large\else
  \ifnum #1<34\Large\else \ifnum #1<41\LARGE\else
     \huge\fi\fi\fi\fi\fi\fi
  \csname #3\endcsname}%
\else
\gdef\SetFigFont#1#2#3{\begingroup
  \count@#1\relax \ifnum 25<\count@\count@25\fi
  \def\x{\endgroup\@setsize\SetFigFont{#2pt}}%
  \expandafter\x
    \csname \romannumeral\the\count@ pt\expandafter\endcsname
    \csname @\romannumeral\the\count@ pt\endcsname
  \csname #3\endcsname}%
\fi
\fi\endgroup
\begin{picture}(6819,5657)(2102,-6993)
\put(3496,-3076){\makebox(0,0)[lb]{\smash{\SetFigFont{10}{12.0}{rm}
$Z_i$
}}}
\put(5515,-4189){\makebox(0,0)[lb]{\smash{\SetFigFont{10}{12.0}{rm}
$\bar\gamma$
}}}
\put(5881,-4921){\makebox(0,0)[lb]{\smash{\SetFigFont{10}{12.0}{rm}
$Z$
}}}
\end{picture}